\documentclass[12pt,leqno]{amsart}  
\usepackage{amsmath,amstext,amsthm,amssymb,amsxtra}
\usepackage[colorlinks,citecolor=red,pagebackref,hypertexnames=false]{hyperref}
\usepackage[backrefs]{amsrefs}
\allowdisplaybreaks
\swapnumbers

\numberwithin{equation}{section}
\theoremstyle{plain} 
\newtheorem{lemma}[equation]{Lemma} 
\newtheorem{proposition}[equation]{Proposition} 
\newtheorem{theorem}[equation]{Theorem} 
\newtheorem{corollary}[equation]{Corollary} 
\newtheorem{conjecture}[equation]{Conjecture}
\newtheorem*{ww}{Wiener--Wintner Theorem}
\newtheorem*{carleson}{Carleson's Theorem}

\theoremstyle{definition}

\theoremstyle{remark}

\newtheorem*{Acknowledgment}{Acknowledgment}

\def\norm#1.#2.{\lVert#1\rVert_{#2}}
\def\Norm#1.#2.{\bigl\lVert#1\bigr\rVert_{#2}}
\def\NOrm#1.#2.{\Bigl\lVert#1\Bigr\rVert_{#2}}
\def\NORm#1.#2.{\biggl\lVert#1\biggr\rVert_{#2}}
\def\NORM#1.#2.{\Bigl\lVert#1\Bigr\rVert_{#2}}

\def\ip#1,#2,{\langle #1,#2\rangle}
\def\Ip#1,#2,{\bigl\langle#1,#2\bigr\rangle}
\def\IP#1,#2,{\Bigl\langle#1,#2\Bigr\rangle}

\def\mid{\,:\,}

\def\abs#1{\lvert#1\rvert}
\def\Abs#1{\bigl\lvert#1\bigr\rvert}
\def\ABs#1{\Bigl\lvert#1\Bigr\rvert}
\def\ABS#1{\biggl\lvert#1\biggr\rvert}

\def\XXint#1#2#3{{\setbox0=\hbox{$#1{#2#3}{\int}$}
     \vcenter{\hbox{$#2#3$}}\kern-.5\wd0}}

\def\name#1 #2{\operatorname{#1}(#2)}
 \def\Tree{{\mathbf T}}
  \def\tree#1{$#1$tree} 

\def\eqdef{\stackrel{\mathrm{def}}{{}={}}}

\begin{document}

\title[Wiener-Wintner for Hilbert Transform] {Wiener-Wintner for Hilbert Transform   }
  
%
%
%
%

\author{Michael Lacey}

\address{Michael Lacey\\
School of Mathematics\\
Georgia Institute of Technology\\
Atlanta,  GA 30332 USA}

\thanks{Research supported in part by the National Science Foundation.  The author 
is a Guggenheim Fellow.}

\email{lacey@math.gatech.edu}

\author{Erin Terwilleger}

\address{
Erin Terwilleger\\
Department of Mathematics, U-3009\\
University of Connecticut\\
Storrs, CT 06269 USA}

\email{terwilleger@math.uconn.edu}

\begin{abstract} We prove the following extension of the Wiener--Wintner Theorem 
and the Carleson Theorem on pointwise convergence of Fourier series: 
For all measure preserving flows $ (X,\mu , T_t)$ and $ f\in L^p (X,\mu)$, 
 there is a 
set $X_f\subset X $ of probability one, so that for all $x\in X_f$ we have 
\begin{equation*}
\lim _{s\downarrow0} \int _{s<\abs t<1/s} \operatorname e ^{i \theta t} f(\operatorname T_tx)\; \frac{dt}t
\qquad \text{exists for all $\theta$.} 
\end{equation*}
The proof is by way of establishing an appropriate oscillation inequality which is itself 
an extension of Carleson's theorem.
\end{abstract}

\maketitle

 \section{ The Main Theorem}

We are concerned with quantitative inequalities related to the pointwise 
convergence of singular integrals that are uniform with respect to modulation. 
To state our results, 
define   dilation and modulation operators by 
 \begin{align} 
 \label{e.dil}
 \operatorname {Dil} ^{(p)}_s f(x)&{}\eqdef{}{s^{-1/p}} f(x/s),\qquad 0<s,p<\infty .
 \\  \nonumber
 \operatorname {Dil} ^{(\infty)}_s f(x)&{}\eqdef{} f(x/s),\qquad 0<s<\infty . 
 \\ 
 \label{e.mod}\operatorname {Mod}_\xi f(x) &{}\eqdef{}\operatorname e ^{ix \xi }f(x), \qquad\xi \in\mathbb R.
 \end{align}

 Let $K $ be a distribution.  The most important example will be $K_H(y) {}\eqdef \zeta (y)
 \frac1y $, where $ \zeta $ is a smooth, symmetric, compactly supported function.  
This is a distribution associated to a truncation of the Hilbert transform kernel.

Our principal concern is the convergence of  terms $(\operatorname {Dil}^{(1)}_s  K)*f(x) $ in a pointwise sense, and in one that 
is, in addition, uniform over all modulations.  To do this, we use the following definition.   
\begin{equation}  \label{e.osc}
\operatorname {Osc}_n (K\,;\, f)^2 {}\eqdef{} \sum_{j=1}^\infty 
\sup _{{k_j\le{} l<l'<k _{j+1 } }} \Abs{ [(\operatorname {Dil}^{(1)}_{2^{l/n} }  K)-(\operatorname {Dil}^{(1)}_
		{ 2^{l'/n} }K)]*f }^2. 
\end{equation}
 This definition depends upon a choice of  an increasing sequence of integers $k_j\in \mathbb Z$, a dependence that 
 we suppress as relevant constants are independent of the choice of $\{k_j\} $.   It also depends upon a choice of 
 positive integer $n $, which we have incorporated into the notation.  This  only permits dilations of the form 
 $2^{l/n} $ for integers $l $.

 \begin{theorem}\label{t.osc}  Fix a smooth, symmetric, compactly supported function $ \zeta $. 
 For integers $n>0 $ and $1<p<\infty $ there is a constant $C_{n,p, \zeta } $ so that we have the inequality 
 \begin{equation}  \label{e.osc<}
 \Norm \sup_N \operatorname {Osc}_n(K_H\,;\, \operatorname {Mod}_N f).p.\le{}C_{n,p,\zeta }\norm f.p..
 \end{equation}
 The inequality holds for all choices of increasing  sequences $\{k_j\mid j\ge1\} $ satisfying $k _{j+1}\ge{} k_j+n$. 
\end{theorem} 

Our primary interest in this theorem is the corollary below, which is a Hilbert transform  counterpart to the well known 
Wiener--Wintner theorem for ergodic averages.  Deriving the corollary below 
is a standard part of the literature, with the roots of the argument going 
back to Calder\'on \cite{MR37:2939}.  The use of an oscillation inequality to establish 
convergence was introduced by Bourgain \cite{MR1019960}.  Also see the papers of
Campbell et al.~\cite{MR1788042}, and Jones et al.~\cite{MR1645330}. 

\begin{corollary}\label{c.osc} For all measure preserving flows $\{ \operatorname T_t \mid t\in \mathbb R\} $
on a probability space $(X,\mu)$ and functions $f\in L ^{p} (\mu)$, there is a 
set $X_f\subset X $ of probability one, so that for all $x\in X_f$ we have 
\begin{equation*}
\lim _{s\downarrow0} \int _{s<\abs t<1/s} \operatorname e ^{i \theta t} f(\operatorname T_tx)\; \frac{dt}t
\qquad \text{exists for all $\theta$.} 
\end{equation*}
\end{corollary}

This is a  common  extension of two classical theorems: Carleson's Theorem \cite {carleson}
on Fourier series with Hunt's extension \cite {MR38:6296}, and the Wiener--Wintner Theorem 
\cite {MR4:15a} on ergodic averages. 

\begin {carleson} We have the inequality 
\begin{equation*}
\NOrm \sup _{N} \ABs{ \int \operatorname {Mod} _{N} f(x-y) \; \frac {dy}y} .p. 
\lesssim \norm f.p.\,, \qquad 1<p<\infty \,. 
\end{equation*}
\end {carleson}

\begin {ww} 
For all measure preserving flows $\{ \operatorname T_t \mid t\in \mathbb R\} $
on a probability space $(X,\mu)$ and functions $f\in L ^{p} (X,\mu)$, there is a 
set $X_f\subset X $ of probability one, so that for all $x\in X_f$ we have 
\begin{gather*}
\lim _{s\to \infty } s ^{-1}\int _{-s} ^{s}
\operatorname e ^{i \theta t} f(\operatorname T_tx)\; {dt}
\qquad \text{exists for all $\theta$}
\\
\lim _{s\to 0 } s ^{-1}\int _{-s} ^{s}
\operatorname e ^{i \theta t} f(\operatorname T_tx)\; {dt}
\qquad \text{exists for all $\theta$.}
\end{gather*}
\end{ww}

The Wiener--Wintner Theorem can been seen as an extension of the 
Birkhoff Ergodic Theorem.  The Carleson Theorem is 
a deep result from the 60's, and since then several proofs have been offered.  An extensive survey and 
bibliography on this subject can be found in \cite {MR2091007}.

The possibility of extending the Wiener--Wintner Theorem to the setting of 
the Hilbert transform was first raised in the paper of Campbell and Petersen 
\cite {MR89i:28008}.
The specific result proved there was  essentially Carleson's 
Theorem on the integers, with a transference to measure preserving systems. 
Part of this was contained in a prior work of M{\'a}t{\'e} \cite {MR39:701}, 
a work that was  overlooked until  much later.  

Assani \cites {MR1821469,MR1995517} proved our Corollary~\ref{c.osc} 
on different classes of dynamical systems.  Indeed, he formulated the concept 
of a \emph{Wiener--Wintner system}.  In this nomenclature, our corollary states 
that all measure preserving systems are Wiener--Wintner systems. 

Our tool to prove convergence in the Hilbert transform setting is the oscillation inequality (\ref{e.osc<}), an idea first employed in 
ergodic theory in the pioneering work of Bourgain on the ergodic theorem 
along arithmetic sequences \cite {MR1019960}.  The use of oscillation has 
subsequently been systematically studied in e.g.~\cites {MR1788042,MR1645330}
and in references therein.

The main goal of the paper is  a proof of Theorem \ref{t.osc}. Clearly, we follow the lines of a proof of Carleson's Theorem. In particular we employ
the Lacey--Thiele approach \cite{laceythielecarleson} and refine one part of it to deduce our main theorem. 
We will also appeal to the `restricted weak type argument' of 
C.~Muscalu, T.~Tao, and C.~Thiele 
\cite{MR2003b:42017}  and L.~Grafakos, T.~Tao, and E.~Terwilleger \cite{gtt}.

\begin {Acknowledgment} The authors have benefited from conversations with 
Jim Campbell, Anthony Quas and Mate Wierdl. Part of this research was completed at the 
Schr\"odinger Institute, Vienna Austria. 
For one of us (ML), discussions with Karl 
Petersen about this question formed our introduction to Carleson's theorem, 
for which we have been indebted to him ever since.
\end {Acknowledgment}

\section{Deduction of  Theorem~\ref{t.osc} }

There are two more technical estimates that we prove. 
Specifically, let $\psi$  be some  Schwartz function which satisfies
\begin{gather}
\label{e.zc0}  0\le \widehat \psi (\xi ) \le{}C_0, 
\\  
\label{e.zc1} 
\widehat \psi \text{\quad is supported in $[-2  ,-\tfrac12 ]$, }
\\
\label{e.zc2} 
\abs {\psi(y) }\le{} C_1\min(\abs y ^{-\nu},\abs y ^\nu ).
\end{gather}
Here, $\nu $ will be a large constant whose exact value we need not specify. 
And we will not have complete freedom in precisely which Schwartz function 
$ \psi $ we can take here.  It should arise in a particular way described in 
the proof of Proposition~\ref{p.osc}, and will be nonzero!  The purpose 
of this section is to describe how a particular result for \emph {any choice 
of non zero  $ \psi $} as above will lead to a proof of our main theorem.

Consider the distribution 
\begin{equation*}
\Psi {}\eqdef{}\sum _{v=1 }^\infty \operatorname {Dil} _{2^{-v} }^{(1) } \psi.
\end{equation*}
We will prove the following two propositions in the next section.

\begin{proposition}\label{p.osc} With the assumptions (\ref{e.zc1})---(\ref{e.zc2}), the inequality (\ref{e.osc<}) holds with 
$n=1 $ and the distribution $K_H $ replaced by $\Psi $. 
\end{proposition}

\begin{proposition}\label{p.singlescale} We have the inequality 
\begin{equation}  \label{e.singlescale}
\NOrm \sup _{N}\Bigl[ \sum_{j=-\infty} ^{\infty} \abs{ \{\operatorname {Dil} ^{(1)}_{2 ^{j} } \psi\}* 
\operatorname {Mod}_N (f) } ^2 \Bigr] ^{1/2} .p. {}\lesssim{} 
\norm f.p. ,\qquad 1<p<\infty. 
\end{equation}
\end{proposition}

 Note that for fixed modulation, (\ref{e.singlescale}) is a Littlewood Paley inequality, making the inequality above 
a ``Carlesonized Littlewood Paley''{} inequality.  Inequalities like this 
have been proved by Prestini and Sj\"olin \cite{MR1781088}. They also follow from the 
method of Lacey and Thiele. 

Both propositions follow from our  Proposition~\ref{p.osc-tiles} of the next section, which is phrased in a language 
conducive to the  methods of Lacey and Thiele \cite{laceythielecarleson}. 
These methods have been applied in a number of variants of Carleson's theorem, 
see e.g. Pramanik and Terwilleger \cite{MR2007237} and Grafakos, Tao and Terwilleger \cite{gtt}.

\smallskip

We turn to the deduction of Theorem~\ref{t.osc}. 
Observe that the two previous propositions immediately prove that when we consider dilations which are 
powers of $2 ^{1/n}$ we have 
\begin{equation*}
\norm \sup_N \operatorname {Osc}_n (\Psi ; \operatorname {Mod}_N f) .p. {}\lesssim{}n\norm f.p.,\qquad n\in \mathbb N,\ 1<p<\infty. 
\end{equation*}
Thus we need not concern ourselves with this feature of Theorem~\ref{t.osc} and Corollary~\ref{c.osc}.  

\medskip 
For a distribution $K$, set 
\begin{equation*}
\norm K.*,p.=\sup _{\norm f.p.=1}\norm \sup_N \operatorname {Osc}_1 (K ; \operatorname {Mod}_N f) .p..
\end{equation*}
Note that since our definition incorporates differences, this is a seminorm on 
distributions $K$. That is, it obeys the triangle inequality (which we use), but can be zero for non zero 
distributions.  In 
particular, for a Dirac point mass $\delta$ we have $\norm \delta.*,p.=0$, and  similarly for the distribution 
$K$ with $\widehat K=\mathbf 1 _{[0,\infty)}$.

Our task is to show that 
$
\norm K_H.*,p.<\infty\,
$, where $ K _{H}(y)=y ^{-1} \zeta (y)$ for some smooth symmetric, compactly supported 
Schwartz function. 
Our Proposition~\ref{p.osc}  is, with this notation, the assertion that $\norm \Psi.*,p.<\infty$.  
The same inequality will hold for a kernel which   can be obtained as a 
convex combination of dilations of $\psi$ and $\Psi$. Thus, set 
\begin{equation*}
 \Psi_0 {}\eqdef{}\int_0^1 \operatorname {Dil}^{(1)}_{2^s} \Psi  \; \frac{ds}s.
\end{equation*}
In this integral, we are careful to integrate against the measure $\frac{ds}s$, which is the Haar measure 
for the positive reals under multiplication, the underlying group for the dilation operators.  In particular, 
it follows that $\Psi_0$ is a distribution whose Fourier transform  is a
nonzero constant on $(-\infty,-1)$ 
and is $0$ on $(-\frac12,\infty)$. 
Thus by Proposition \ref{p.osc}, we clearly have $\norm \Psi_0.*,p.<\infty$.

Now we will show that $\norm D_0.*,p.<\infty$ for the distribution 
\begin{equation*}
D_0(y)=y ^{-1}\zeta (y)-c(\Psi_0(y)-\overline{ \Psi_0(y)}),
\end{equation*}
where we choose the complex constant $c$ so that $\lim _{ \xi \rightarrow \infty} \widehat D_0(\xi)=0.$
In fact, it is a well known elementary fact that for $c=i\pi$, 
\begin{equation} \label{e.decay}
\int \zeta (y) \operatorname e ^{i \xi y}\frac {dy}y=c+O(\abs{ \xi} ^{-1}). 
\end{equation}  
We will decompose the distribution $D_0$ into a sum which
can be treated with Proposition~\ref{p.singlescale}.  
Then using that $\norm \Psi_0.*,p.<\infty$ and $\norm D_0.*,p.<\infty$, 
we obtain the desired inequality for $K(y)=y ^{-1}\zeta (y)$.

Choose $\chi $ to be a smooth function supported on $\tfrac12\le \abs{ \xi }\le2$
so that
\begin{equation*}
\sum_{k=-\infty} ^{ \infty } \operatorname {Dil} ^{ \infty } _{ 2^{-k}}\chi=\mathbf 1 _{\mathbb R-\{0\}}, 
\end{equation*}
and set $\widehat \Delta_k=\widehat D_0 \operatorname {Dil} ^{ \infty } _{ 2^{-k}}\chi$.  
The following lemma finishes the proof of Theorem~\ref{t.osc}.

\begin{lemma}\label{l.Delta}  
We have 
\begin{equation*}
\norm \Delta_k.*,p. {}\lesssim{} 2 ^{ -\abs{ k}},\qquad k\in \mathbb Z.
\end{equation*}
\end{lemma}

\begin{proof}
We will verify that 
\begin{gather}\label{e.Delta0}
\norm \widehat \Delta_k .\infty. {}\lesssim{}2 ^{ -|k|},\qquad k\in \mathbb Z\,,
\\ \label{e.Delta1}
\text{ $\widehat \Delta_k$ is supported on $ 2 ^{-k-1}\le \abs{ \xi }\le 2 ^{-k+1}$, }
\\ \label{e.Delta2}
\abs{ \Delta_k (y)} {}\lesssim{} 2 ^{ - k-\abs{ k}}( 1+2 ^{- k }\abs y) ^{ -\nu}\,,
\qquad k\in \mathbb Z,\ y\in \mathbb R\,,
\end{gather}
with implied constants independent of $ k\in \mathbb Z$ and $\nu$ the 
large, unspecified constant that appears in (\ref{e.zc2}). 
With decay in $ \abs{ k}$ in both (\ref{e.Delta0}) and (\ref{e.Delta2}), the lemma then follows from 
a trivial change of scale and 
from Proposition~\ref{p.singlescale}.

Let us recall the trivial estimate which follows from the symmetry of $ \zeta $,  
\begin{equation}\label{e.KH1}
\abs{ \widehat{K_H}(\xi )}=\ABs{\int  \zeta  (y)\frac{ \operatorname e ^{i \xi y } }y \; dy}  \lesssim \abs{ \xi }. 
\end{equation} 
In addition we have the estimate below,  applied for $ \abs{ \xi }\le1$
\begin{equation}\label{e.KH2}
\ABs{ \tfrac {d^w}{d \xi ^w}\widehat{K_H}(\xi )}
=\ABs{\int \zeta (y)y ^{w-1}\operatorname e ^{i \xi y } \; dy} 
\lesssim 
\begin{cases}
\abs{ \xi } & \text{$ w$ even }
\\
1 & \text{$ w$ odd. }
\end{cases}
\end{equation}
Whereas for $ \abs{ \xi }\ge1$, we have 
\begin{equation}\label{e.KH3}
\ABs{ \tfrac {d^w}{d \xi ^w}\widehat{K_H}(\xi )} \lesssim \abs{ \xi } ^{-\nu}, \qquad \abs{ \xi }>1\,,\ 
0<w\le \nu \,.
\end{equation}
That is, we have very rapid decay in a large number of derivatives. 
\smallskip

Now, (\ref{e.Delta1}) is true by definition of $ \Delta _k$.   To see (\ref{e.Delta0}) 
for $ k\ge2$, note that this is only determined by the Fourier transform of $ K_H$ since $\widehat{\Psi_0}$ and $\widehat{\overline{\Psi_0}}$ are zero. The result easily follows by the inequality (\ref{e.KH1}) and property (\ref{e.Delta1}).  For $ k\le2$, the inequality follows from 
the construction of $ D_0$, and in particular the property in (\ref{e.decay}).

We turn to the last condition, (\ref{e.Delta2}). 
It is well known that decay of order $\nu$ in spatial variables is implied by 
differentiability of a function in frequency variables.  
Observe that 
\begin{align*}
\bigl( y ^{\nu} \Delta _k(y)\bigr)\widehat {\;} (\xi ) &=i ^{-\nu} \tfrac {d ^{\nu}} {d \xi ^{\nu}} 
\widehat {\Delta _{k}}(\xi )
\\
&=
i ^{-\nu} \tfrac {d ^{\nu}} {d \xi ^{\nu}} \operatorname {Dil} _{2 ^{-k}} ^{(\infty )} 
\chi (\xi ) \widehat {K_H}(\xi ).  
\end{align*}
Hence, 
\begin{equation*}
\Abs{\bigl( y ^{\nu} \Delta _k(y)\bigr)\widehat {\;} (\xi )}\le
\sum _{w=0} ^{\nu} 2 ^{k w} \ABs{ \tfrac {d ^{\nu-w}} 
{d \xi ^{\nu-w}} \sup _{2 ^{-k-1}\le \abs{ \xi }\le 2 ^{-k+1}}\widehat {K_H}(\xi )  }.
\end{equation*}

For $ k\ge1$, this sum is dominated by the last two terms.  To control them, use (\ref{e.KH2}),
supplying the estimate  
$ \lesssim  2 ^{(\nu-1)k}$.  This is better by a factor of $ 2 ^{-k}$ than the trivial estimate, so 
that Fourier inversion proves (\ref{e.Delta2}) in this case. 

The case of $ k\le0$ is easier, due to the rapid decay in (\ref{e.KH3}). 
\end{proof}

\section{Decomposition and Main Proposition}

We state the definitions needed for the main proposition and conclude this section with the 
argument of how this proposition proves the results of the previous section, namely Proposition~\ref{p.osc} and Proposition~\ref{p.singlescale}. 

In addition to the modulation and dilation operators in
(\ref{e.dil}) and (\ref{e.mod}), we need translation operators 
\begin{equation}  \label{e.trans}
\operatorname {Tran}_y f(x) {}\eqdef{}f(x-y),\qquad y\in \mathbb R . 
\end{equation}
We set $\mathcal D $ to be the dyadic grid and say that $I\times\omega \in\mathcal D\times\mathcal D $ is a {\em tile} iff 
$\abs{ \omega}\cdot\abs I=1 $.  Let $\mathcal T $ denote the set of all tiles. 

We think of $\omega $ as a frequency interval and $I$ as a spatial interval;
our definition of 
a tile is a reflection of the  uncertainty principle for the Fourier transform.   
We  will plot frequency intervals in the vertical direction.  Each dyadic interval $\omega $ is a union of 
two dyadic intervals of half the length of $\omega $.  We call them $\omega_\pm $ and view $\omega_+ $ as above $\omega_- $. 

We take a fixed Schwartz function $\varphi $ with frequency support in the interval $[-1/\nu,1/\nu] $.  
For a tile $s=I_s\times \omega_s $, define 
\begin{equation}  \label{e.zvf}
\varphi_s {}\eqdef{}\operatorname {Mod}_{c(\omega _{s-}) } \operatorname {Tran}_{c(I_s)} \operatorname {Dil}_{\abs {I_s} }^2 
\varphi.  
\end{equation}
Here, $c(J)$ is the center of the interval $J $, and $\omega_{s- } $ is the lower half of the interval $\omega_s $.  Thus, 
this function is localized to be supported in the time frequency plane close to the rectangle $I_s\times \omega _{s-} $. 

There are companion functions which depend
on different choices of certain measurable functions.  These functions should be 
thought of as those choices of modulation and indices that will achieve, up to a constant multiple, the supremums 
in the oscillation function.  To linearize the modulation, let 
\begin{equation*}
\text{$N \,:\,\mathbb R\longrightarrow \mathbb R $  be a measurable function (a modulation parameter)}.   
\end{equation*}
We define another function  related to the rectangle $I_s\times \omega _{s+} $ which tells us when the linearized modulation parameter is at a certain frequency.  Let
\begin{equation} \label{e.zf}
	\phi_s(x) {}\eqdef{}
	 \mathbf 1 _{\omega_{s+}}(N(x)) \varphi_s (x).
\end{equation}
Now define a tile variant of the oscillation operator by 
\begin{equation}  \label{t.osc-tiles}
\operatorname {Tile-osc}(f) {}\eqdef{}
\Bigl[\sum_{j=1}^\infty 
\sup _{{k_j\le{} l<l'<k _{j+1 } }} \ABs{ \sum _{\substack{s\in\mathcal T \\ 2^l \le{}\abs {I_s }\le 2 ^{l ' } } } 
\ip f,\varphi_s , \phi_s }^2 \Bigr]^{1/2 }.
\end{equation}
Here, an increasing sequence of integers 
$\{k_j\mid j\ge1\} $ are specified in advance. 
We make the definition for clarity's sake, as we will not explicitly work with it.  Rather we 
prefer to fully linearize this maximal operator.  This requires the additional choices of 
functions 
\begin{gather}
\alpha _{j } \,:\,\mathbb R\longrightarrow \mathbb R ,\qquad \sum _{j=1}^\infty \abs{ \alpha _{j }(x) }^2\le{}1, 
\quad \text{for all $x $ }, 
\\
\ell_{j-},\ell _{j+ }\,:\,\mathbb R\longrightarrow \mathbb Z,\quad k_j\le{}\ell _{j- }<\ell _{j+ }<k_{j+1}.
\end{gather}
And we set
\begin{gather} 
\label{e.f}
	F_{s,j} {}\eqdef{}\{ x\mid 2 ^{\ell _{{j-} }(x) }\le{} \abs {I_s } 
	<{}2 ^{\ell _{{j+} }(x) } \}, 
\\
	f_{s,j}(x){}\eqdef{} \mathbf 1 _{F_{s,j}}(x) \alpha_j(x) \phi_s (x) .
\end{gather}
The sequences of functions $\ell_{j\pm}$ are selecting the level at which the maximal difference occurs. 
The $\alpha_j $ are chosen to realize the $\ell^2 $ norm in the definition of oscillation.  
We make all of these choices in order to linearize the oscillation operator.

Our main proposition is 
\begin{proposition}\label{p.osc-tiles}
For all choices of $N(x) $ and  increasing sequence of integers $\{ k_j \} $, the operator $\operatorname {Tile-osc } $ 
extends to a bounded sub linear operator on $L^p $, $1<p<\infty $.   In particular, for sets $G,H\subset \mathbb R $ 
of finite measure, we have 
\begin{equation}  \label{e.osc-tiles}
\sum_{\substack{s\in \mathcal T\\ j\in \mathbb N } }
\abs{\ip \mathbf 1 _{G},\varphi_{s} , \ip \mathbf 1 _{H}, f_{s,j} , } {}\lesssim{}
\min(\abs G,\abs H) \big( 1+\big|\log \tfrac{\abs G}{\abs H}\big| \big) .
\end{equation}

\end{proposition}

Note that the  inequality above implies that 
\begin{equation*}
 \abs{\ip \name {Tile-osc} {\mathbf 1 _{G}}, \mathbf 1 _{H} , } {}\lesssim{} \abs{G} ^{1/p}\abs H  ^{1-1/p},\qquad 1<p<\infty.
\end{equation*}
That is, we have the restricted weak type inequality for all $1<p<\infty $.\footnote{In fact, the estimate (\ref{e.osc-tiles})
gives a favorable upper bound on the behavior of the constant with respect to $p$, namely that they are 
no more than $\max(p,\tfrac p{p-1})$. See \cite{gtt}} Hence, an interpolation argument will give us the estimate 
\begin{equation}  \label{e.osc-tiles-Lp}
\norm \name {Tile-osc} f.p. {}\lesssim{}\norm f.p.,\qquad 1<p<\infty. 
\end{equation}

 \subsection*{The Deduction of Proposition~\ref{p.osc} and Proposition~\ref{p.singlescale}.} 
 
 For $\xi \in \mathbb R$ and $\ell\in \mathbb Z$, consider the operators 
 \begin{equation*}
 \operatorname A _{\xi,l} f {}\eqdef{}\sum _{\abs{I_s}=2 ^{l}} \mathbf 1 _{\xi\in \omega_{s+}}
 \ip f,\varphi_s,\, \varphi_s. 
 \end{equation*}
 The tile oscillation operator is built up from these operators.  Observe that these 
operators enjoy the properties 
\begin{align}\label{trans}
\operatorname A _{\xi,l} \operatorname {Trans}_{n2^l}  & 
	{}={} \operatorname {Trans}_{n2^l} \operatorname A _{\xi,l}, \qquad n \in \mathbb Z ,
\\\label{dil}
\operatorname A _{\xi,l} \operatorname {Dil}^{(2)}_{2 ^{-l ' }}& 
	{}={} \operatorname {Dil}^{(2)}_{2 ^{-l ' }} 
	\operatorname A _{\xi 2 ^{-l'},l+l'}, \qquad l'\in\mathbb Z,
\\ \label{mod}
\operatorname A _{\xi,l} \operatorname {Mod}_{-\theta} & 
 	{}={} \operatorname {Mod}_{-\theta} \operatorname A _{\xi+\theta,l},
	\qquad \theta \in \mathbb R.
\end{align}
Notice that these conditions tell us that the operators $\operatorname A _{\xi,l}$ have a near
translation invariance, a certain modulation invariance, and are related to each other through
dilations. In addition, these operators are bounded on $L^2$ uniformly in $\xi$ and $l$, a fact
well represented in the literature.

We will now define 
\begin{align*}
\operatorname B _{\xi,l} f&{}\eqdef{}
\lim_{\substack {K \to \infty \\ L\to \infty}}
\frac{1}{4KL}\int_{-K}^{K}\int_{-L} ^{L}
\operatorname{Mod}_{-\theta}\operatorname{Trans}_{-y}
\operatorname A _{(\xi+\theta),l}
(
\operatorname {Trans}_{y}	
\operatorname {Mod}_{\theta}f)\, dy d\theta\,. 
\end{align*}
By periodicity of the integrand in $y$ and $\theta$, for all Schwartz functions $ f$, 
the averages on the right hand 
side converge pointwise to  $\operatorname B _{\xi,l}f(x)$ as $K,L \to \infty$.  

Let us make some observations about the operators $\operatorname B _{\xi,l}$.  First,
 (\ref{trans}) and periodicity of the integrand in $y$ imply $\operatorname B _{\xi,l}$
 commutes with translations. Second, it is a bounded, positive, semidefinite operator, 
 as is easy to see. 
  Hence, it is given by convolution.  Indeed, (\ref{mod}) implies that 
  \begin{equation*}
\operatorname B _{\xi ,l} f= \operatorname {Mod} _{\xi } \beta _{l} \ast (
\operatorname {Mod} _{-\xi } f)\,.
\end{equation*}
 for a function $ \beta _{l}$ that we turn to next.  
 The equality (\ref{dil}) implies that $ \beta _{l}= \operatorname {Dil} _{2 ^{l}} ^{(1)} \beta
 $ where $\beta$ is given such that $ \beta _{0}$ is 
 a smooth Schwartz function satisfying the conditions (\ref{e.zc0})---(\ref{e.zc2}), a routine exercise to verify.

 Assuming Proposition~\ref{p.osc-tiles},  it follows that we can conclude 
 Proposition~\ref{p.osc} and Proposition~\ref{p.singlescale} for nonzero functions 
 $ \psi =\beta_0 $. Our proof is complete.

\section{Main Lemmas} 

To prove (\ref{e.osc-tiles}), we split the sum over $s\in\mathcal T$ into the sum over $s$  such that $I_s \subset \{\operatorname M \mathbf 1 _{G}>\lambda \}$ and the sum over $s$ such that $I_s \not\subset \{\operatorname M \mathbf 1 _{G}>\lambda \}$. The former sum can be taken care of by an argument of M. Lacey and C. Thiele \cite{lacey-thiele2} which also appears, slightly modified, in the paper of L. Grafakos, T. Tao, and E. Terwilleger \cite{gtt}. Thus we restrict our attention to the tiles $s$ where $I_s \not\subset \{\operatorname M \mathbf 1 _{G}>\lambda \}$.

We begin with some concepts needed to phrase the proof.  There is a natural partial order on tiles. We say that $s<s'$ iff $\omega_{s}\supset 
\omega_{s'}$ and $I_s\subset I_{s'}$.  Note that the time variable of $s$ is localized to that of $s'$, and the frequency 
variable of $s$ is 
similarly localized, up to the variability allowed by the uncertainty principle. Note that two tiles are incomparable 
with respect to the 
`$<$' partial order iff the tiles, as rectangles in the time frequency plane, do not intersect.   A ``maximal tile''{} will be
one that is maximal 
with respect to this partial order.   

Let $\mathcal S$ denote an arbitrary set of tiles.  We call a set of tiles $\Tree \subset\mathcal S$ a {\em tree} if there is a tile $I_\Tree\times \omega_{\Tree}$, called the
{\em top of the tree,} 
such that for all $s\in\Tree$, $s<I_\Tree\times \omega_{\Tree}$.  We note that the top is not uniquely defined.  An important point is that a tree top specifies a location in time variable for the tiles in the tree, namely inside $I_\Tree$, and localizes the  frequency variables,  identifying $\omega_{\Tree}$ as 
a nominal origin.

We say that {\em the count of $\mathcal S $ is at most $A $ } iff 
 $\mathcal S=\bigcup_{\Tree\subset\mathcal S}\Tree$, where each $\Tree\subset\mathcal T$ is a tree which is maximal with respect to inclusion and 
\begin{equation*}
\name Count {\mathcal S} {}\eqdef{}\sum_{\Tree\subset\mathcal S}\abs {I_\Tree }\le{} A .
\end{equation*}

Fix $\chi(x)=(1+\abs x)^{-\nu}$, where $\nu$ is, as before, a large constant whose exact value is unimportant to us.  Define  
\begin{gather} \label{e.chiI} 
\chi_I:=\operatorname  {Trans}_ {c(I)} \operatorname { Dil} _{\abs I}^{(1)}  \chi,
\\
\label{e.dense} 
\name dense  s:=\sup_{s<s'}\int_{ N^{-1}(\omega_{s'})\cap H} \chi_{I_{s'}}\; dx,
\\\nonumber 
\name dense  {\mathcal S}:=\sup_{s\in\mathcal S}\name dense  s,\qquad \mathcal S\subset\mathcal T.
\end{gather}

The first and most natural definition of a ``density'' of a tile, would be ${\abs {I_s}}^{-1}\abs{N^{-1}(\omega_{s+})\cap I_s}$.  However 
$\varphi$ is supported on the whole real line, although it does decay faster than the inverse of any  polynomial. 
We refer to this as a ``Schwartz tails problem.''
  The definition of density as  $\int_{N^{-1}(\omega_{s})} \chi_{I_s}\;dx $, as it turns out, is still not adequate. 
 That we should take the supremum over  $s<s'$ only becomes evident in the proof of the ``Tree Lemma'' below.  
 
 The ``Density Lemma'' is 
 

 \begin{lemma}\label{l.dense} Any subset $\mathcal S\subset\mathcal T$ is a union of $\mathcal S_{\text{\rm heavy}}$ and $\mathcal S_{\text{\rm light}}$ for which 
 \begin{equation*}
 \name dense  {\mathcal S_{\text{\rm light}}}<\tfrac12\name dense { \mathcal S} ,
 \end{equation*}
 and the collection $\mathcal S_{\text { \rm heavy } }$  satisfies 
 \begin{equation}\label{e.densecount}
 \operatorname{Count} (\mathcal S_{ \text{\rm heavy }}){}\lesssim{}{}{\name dense  {\mathcal S}}^{-1} \abs H .
 \end{equation}
 \end{lemma}

  What is significant is that this relatively simple lemma admits a non-trivial variant  intimately
  linked to the tree structure and orthogonality.  We should refine the notion of a tree. Call a tree $\Tree$ with top $I_\Tree\times \omega_{\Tree}$
 a \tree\pm{} iff for each $s\in\Tree$, aside from the top, 
  $I_\Tree\times \omega_{\Tree}\cap I_s\times \omega_{s\pm}$ is not empty.  Any tree is a union of a 
  \tree+{} and a \tree-{}.  If $\Tree$ is a \tree+{}, 
  observe that the rectangles $\{I_s\times \omega_{s-}\mid s\in\Tree\}$ are disjoint.  
  We see that 
  \begin{equation*}
  \sum_{s\in\Tree}\abs{\ip f,\varphi_s, }^2{}\lesssim{}{}\lVert f\rVert_2^2.
  \end{equation*}
  This motivates the  definition 
  \begin{equation} \label{e.size}
  \name size  {\mathcal S}:=\sup\{ \abs{I_\Tree}^{-1/2} 
  \sum_{s\in\Tree}\abs{\ip f,\varphi_s, }^2\mid \Tree\subset\mathcal S,\ \text{$\Tree$ is a \tree+{}}\}.
  \end{equation}

  The ``Size Lemma'' is 
  
 
 \begin{lemma}\label{l.size} Assume that $f=\mathbf 1 _{ G} $.  
 Any subset $\mathcal S\subset\mathcal T$ is a union of $\mathcal S_{\text{\rm big}}$ and $\mathcal S_{\text{\rm small}}$ for which 
 \begin{equation*}
 \name size  {\mathcal S_{\text{\rm small}}}<\tfrac12\name size  {\mathcal S} ,
 \end{equation*}
 and the collection $\mathcal S_{\text{\rm big}}$ satisfies 
 \begin{equation} \label{e.sizecount} 
  \operatorname{Count} ({\mathcal S_{\text{\rm big}}}){}\lesssim{}{}{\name size  {\mathcal S}}^{-2} \abs G .
 \end{equation}
 \end{lemma}

 
 Concerning the quantity size, we need an additional piece of information about it.  Recall that $\operatorname M$ is the Hardy Littlewood maximal function. 
 
 \begin{lemma}\label{l.size-upper} Let  $0< \lambda<1 $, and suppose that $\mathcal S $ is the set of tiles with 
 \begin{equation*}
 I_s \not\subset \{\operatorname M \mathbf 1 _{G}>\lambda \},\qquad s\in\mathcal S. 
 \end{equation*}
 Then it is the case that $\name size {\mathcal S} {}\lesssim{}\lambda $.
 \end{lemma} 
 
 This fact, a delicate consequence of the Calderon--Zygmund decomposition, will not be proved in this paper. 
 It, like the Size Lemma and the Density Lemma, is already well represented in the literature.  
   See, for example,  \cite{gtt}.  For proofs of the Density and Size Lemmas, we 
   refer the reader to \cite{laceythielecarleson}.  The survey \cite {MR2091007} 
   is also suggested.


For a set of tiles $\mathcal S $, set 
\begin{equation*}
\operatorname {Sum}(\mathcal S) {}\eqdef{}\sum_{j\in \mathbb N}\sum_{s\in\mathcal S} \abs{ \ip \mathbf 1 _{G} , \varphi_s , \ip \mathbf 1 _{H}, f_{s,j} ,}
 \end{equation*}
 Our final lemma relates trees, density and size. It is the ``Tree Lemma.'' 
 
 \begin{lemma}\label{l.tree} For any tree $\Tree$ 
 \begin{equation}\label{e.tree}
 \name Sum {\Tree}{}\lesssim{}{}\name size  \Tree \name dense \Tree\abs{I_\Tree} .
 \end{equation}
 \end{lemma}
 
 Of course for any set of tiles $\mathcal S $, we would then have 
 \begin{equation*}
  \name Sum {\mathcal S}{}\lesssim{}{} \sum_{\Tree\subset\mathcal S}\name size  \Tree \name dense \Tree |I_{\Tree}| .
 \end{equation*}
 Thus, we should inductively apply 
 Lemma~\ref{l.dense} and Lemma~\ref{l.size} so that the `Count' estimates are essentially equal.  
 The formal proof of Proposition~\ref{p.osc-tiles}, which is 
much as it appears in Lacey and Thiele \cite{laceythielecarleson} with the adaptation to a restricted weak type inequality as seen in \cite{gtt}, is left as an exercise for the reader.

 \section{Proof of Lemma~\ref{l.tree}}   
 
The tree lemma, with its adaptation to the setting of oscillation,  is the 
primary new step in this paper.

We begin with some remarks about oscillation operators, and a particular form of the same that we shall use 
at a critical point of this proof.   Let $\zeta$ be a smooth function with Fourier transform supported in
$[-1-\varepsilon,1+\varepsilon]$  for a fixed, small, positive $\varepsilon$ and equal to 1 on $[-1,1]$. 
Set 
\begin{align*}
\operatorname {Osc} (f)^2 & {}\eqdef{}  \sum_{j=1}^\infty \sup _{2 ^{k_j}\le{}\abs I \le{}\abs {I'} \le{}  2 ^{k_{j+1}
}  } \abs{  \operatorname {Dil}^{(1)}_{\abs I}\zeta*f-{}
\operatorname {Dil}^{(1)}_{\abs{ I' }}\zeta*f }^2.
\end{align*}
It is  known that this is bounded on $L^2$, and in this situation we will give an elementary proof of this fact below.

We shall have recourse to  not only this bound, but a particular refinement.
Let $\mathcal J$ be a partition of $\mathbb R$ 
into dyadic intervals.  To each $J\in\mathcal J$, associate a subset 
$E(J)\subset J$ with $\abs{E(J)}\le{}\delta\abs J$, 
where $0<\delta<1$ is fixed. 
Consider 
\begin{equation}  \label{e.oscWithDenseDef}
\operatorname {Osc}_\delta (f)^2 {}\eqdef{} 
\sum _{J\in \mathcal J}  \mathbf 1 _{E(J)}
\sum_{j=1}^\infty \sup _{\substack{J\subset I \subset I'\\ \substack{2 ^{k_j}\le{}\abs I \le{}\abs {I'} \le{}  2 ^{k_{j+1}
}  }} } \abs{ \mathbf 1 _{I}\ip \zeta_I,f,-\mathbf 1 _{I'}\ip \zeta_{I'} ,f, }^2
\end{equation}
We estimate the norm of this operator.

\begin{lemma}\label{l.oscWithDense} We have the estimate 
\begin{equation} \label{e.oscWithDense}
\norm \operatorname {Osc}_\delta (f).2. \lesssim \sqrt \delta \norm f.2.
\end{equation}
for all $f\in L^2$.
\end{lemma}

\begin{proof}
Let us begin with a proof that $ \norm \operatorname {Osc}.2\to2. \lesssim 1$.  
That is, we do not have the additional information about the partition $ \mathcal J$, 
and sets $ E (J)$ for $ J\in \mathcal J$. 
For a sequence of increasing integers $ k_j$ and function $ f\in L^2$, set  
\begin{equation*}
 \widehat { f_j}= \mathbf 1 _{2 ^{-k_{j+1}-1}\le \abs{  \xi } \le 2 ^{-k _{j}+1} } 
\widehat f  
\end{equation*}
Then, we certainly have $ \sum _{j\in\mathbb N} \norm f_j.2.^2\le{} 3\norm f.2.^2$.
Moreover,
due to our assumption about the function $ \zeta $, 
\begin{equation*}
\sup _{2 ^{k_j}\le{}\abs I \le{}\abs {I'} \le{}  2 ^{k_{j+1}
}  } \abs{  \operatorname {Dil}^{(1)}_{\abs I}\zeta*f}
\le
\operatorname M f _{j-1}+\operatorname M f _{j}+\operatorname M f _{j+1}
\end{equation*}
where $ \operatorname M$ is the usual maximal function.  Thus, by the boundedness 
of the maximal function on $ L^2$ we have 
\begin{equation*}
\norm \operatorname {Osc}(f).2.^2
\le
3\sum _{j} \norm \operatorname M f_j.2.^2 \lesssim \norm f.2.^2\,.
\end{equation*}

It is hardly surprising that the proof above appeals to the boundedness of the maximal function, 
as the estimate on the oscillation operator implies that for the maximal function.  
Likewise, our lemma implies a bound for a certain variant of the maximal function. 
As it turns out, we need this variant in the course of the proof. 

Define 
\begin{equation*}
\operatorname M _{\delta } f(x) \eqdef \sum _{J\in \mathcal J} \mathbf 1 _{E(J)}(x) \sup _{J\subset I} 
\ip \abs  f ,\chi _I, 
\end{equation*}
where $ \chi _I$ is defined as in (\ref{e.chiI}).  Then  the estimate we claim is
$ \norm \operatorname M _{\delta }.2. \lesssim \sqrt \delta $.  

Indeed, for any  point $ x\in E(J)$, we have the inequality 
\begin{equation*}
 \operatorname M _{\delta }f(x) \lesssim 
\inf _{y\in J}\operatorname M f(y).
\end{equation*}
where $ \operatorname M$ is the usual maximal function. 
Therefore, we can estimate 
\begin{align*}
\norm \operatorname M _{\delta }f.2.^2&= 
	\sum _{J\in \mathcal T} \int _{E(J) }\operatorname M _{\delta }f(x)^2\; dx
\\
& \lesssim  \sum _{J\in \mathcal T} \abs{ E(J)} \inf _{y\in J}\operatorname M f(y) ^2 
\\
&\le{} \delta \int \operatorname M f(x)^2\; dx \,.
\end{align*}
This proves our claim.

\smallskip 

To conclude the proof, we can estimate 
\begin{align*}
\int_{\cup_{J}E(J)} {\operatorname {Osc}_\delta (f)}(x)^2\; dx\lesssim {}&
\sum_{j}\norm \operatorname M _{\delta } f_j.2.^2 
\\ 
\lesssim & \delta \sum_{j}\norm f_j.2.^2
 \\{}\lesssim{}&\delta\lVert f\rVert_2^2.
 \end{align*}
Our proof is complete.

\end{proof}

We begin the main line of the argument. Let $\delta=\operatorname {dense} (\Tree) $,
 and $\sigma=\operatorname {size} (\Tree) $.  By a modification of the 
 functions $ \alpha _j(x)$ by a choice of signs, we can assume the identity 
 \begin{equation*}
\sum_{j\in \mathbb N} \sum_{s\in\Tree}\abs{\ip \mathbf 1 _{G},\varphi_s ,\ip f_{s,j},{\mathbf 1}_ H ,}={} 
 	\int_H\sum_{j\in \mathbb N}{\sum_{s\in\Tree}  
	\ip \mathbf 1 _{G},\varphi_s, f_{s,j}(x) }\; dx.
\end{equation*}

 As we have no particular control on the set $H$, we will need the following partition of the real line 
 induced by the tree $\Tree$.  
 Let $\mathcal J$ be the partition of $\mathbb R$ consisting of the maximal dyadic intervals 
 $J$ such that $3J$ does not contain any $I_s$ for $s\in\Tree$.  
 It is helpful to observe that for such $J$, if $\abs{J}\le{}\abs{I_\Tree}$, then $J\subset 3 I_\Tree$, and if $\abs J\ge\abs {I_\Tree}$,
 then 
 $\text{dist}(J,I_\Tree)\gtrsim\abs J$. 
 The integral above is at most the sum of the two terms below.
 \begin{gather} \label{e.firstsum}
\sum_{j\in \mathbb N} \sum_{J\in\mathcal J}
\sum_{\substack{s\in\Tree\\ \abs{I_s}\le2\abs J}}\abs{\ip {\mathbf 1}_G, \varphi_s,} 
\int_{J\cap H}\abs{f_{s,j}(x)}\; dx
 \\
 \label{e.secondsum}
\sum_{j\in \mathbb N} \sum_{J\in\mathcal J}\int_{J\cap H}
\ABs{ \sum_{\substack{s\in\Tree\\ \abs{I_s}>2\abs J}}
\ip {\mathbf 1}_G, \varphi_s,  f_{s,j}(x) }\; dx 
 \end{gather}
 Notice that for the second sum to be non--zero, we must have $J\subset 3I_\Tree$.

 The first term (\ref{e.firstsum}) is controlled by an appeal to the ``Schwartz tails.''  Fix an integer $n\ge-1$, and only consider those $s\in\Tree$ for which $\abs{I_s}=2^{-n}\abs J$. 
 Recalling that $f_{s,j}(x){}={} \mathbf 1 _{F_{s,j}}(x) \alpha_j(x) \mathbf 1 _{\omega_{s+}}(N(x)) \varphi_s (x)$, we see that 
\begin{align*}
\sum_{j\in\mathbb N}\sum_{\substack{s\in\Tree\\ \abs{I_s}=2^{-n}\abs J}}
\abs{\ip \mathbf 1 _{G},\varphi_s,} \int_{J\cap H}\abs{f_{s,j}(x)}\; dx\lesssim{}&\sum_{\substack{s\in\Tree\\\abs{I_s}=2^{-n}\abs J}}\sigma\delta (\abs{I_s}^{-1}\text{dist}(I_s,J))^{-10} \abs{I_s}\\
\lesssim{}&\sigma\delta 2^{-n} \text{min}\bigl(\abs{J},\abs J(\abs{I_\Tree}^{-1}\text{dist}(J,I_\Tree))^{-5}\bigr),
\end{align*}
Observe that for each $ s$ above, only one value of $ j$ contributes to the 
left hand sum.  In addition, we have used the 
the fact that there are only a bounded number of tiles 
$s$ for which $\abs{I_s}^{-1}\text{dist}(I_s,J)$ is essentially constant.  
In addition, for the case  $|J|\le|I_{\Tree}|$, we used that
the distance from $I_s$ to $J$ is at least ${}\gtrsim\abs J$.
In the case $|J|>|I_{\Tree}|$, use $\abs{I_s}^{-1}\text{dist}(I_s,J)\ge\abs{I_\Tree}^{-1}\text{dist}(J,I_\Tree)$.
The estimate above can then be summed over $n\ge-1$ and $J\in\mathcal J$ to bound (\ref{e.firstsum}) by ${}{}\lesssim{}\sigma\delta\abs{I_\Tree}$, as required.

Now we turn to the   control of (\ref{e.secondsum}).  
The integral in this quantity is supported in the set  
\begin{equation} \label{e.GJ} 
E(J)=J\cap\bigcup_{\substack{s\in\Tree\\ \abs{I_s}>2\abs J}} \big(N^{-1}(\omega_{s+})\cap H\big).
\end{equation}
Then the critical observation is that $\abs{E(J)}{}\lesssim{}\delta\abs J$. To see this, let $J'$ be the next larger dyadic interval that contains 
$J$.  Then $3J'$ must contain some $I_{s'}$ for $s'\in\Tree$.  Then there exists a tile $s''$ with $I_{s'}\subset I_{s''}\subset I_{\Tree}$ such that
$\abs{I_{s''}}$ is $2\abs J$ or $4\abs J$, and $\omega_\Tree\subset \omega_{s''}\subset\omega_{s'}$.  Then, $s'<s''$, and by the definition of density, 
\begin{equation*}
\frac{|J\cap H\cap N^{-1}(\omega_{s''})|}{|I_{s''}|}{}\le{} \int_{H\cap N^{-1}(\omega_{s''})} \chi_{I_{s''}} \; dx\le\delta
\end{equation*}
But, for each $s$ as in (\ref{e.GJ}), we have $\omega_s\subset \omega_{s''}$, so that
 $E(J)\subset N^{-1}(\omega_{s''})$. Our claim follows.

 Suppose that $\Tree$ is a $-$tree.  This means that the tiles $\{I_s\times \omega_{s+}\mid s\in\Tree\}$ are disjoint  
 and thus the functions $f_{s,j}$ are disjointly supported.  In particular, the oscillation that 
 arises from such functions is trivially bounded by their $\ell  ^{\infty}$ norm.  
 Then the bound for (\ref{e.secondsum}) is no more than 
 \begin{equation*}
 \sum_{j\in\mathbb N}\abs{E(J)} \NOrm \sum_{\substack{s\in\Tree\\ \abs{I_s}>2\abs J}} \abs{\langle {\mathbf 1}_G, \varphi_s\rangle f_{s,j}} .\infty.{}\lesssim{}\delta\sigma\abs{J}.
 \end{equation*}
 This is summed over $J\subset 3I_\Tree$ to get the desired bound.

 Suppose that $\Tree$ is a $+$tree.  This is the interesting case.
 At this point, we will appeal to the norm bound for oscillation, (\ref{e.oscWithDense}), 
 applied to the function 
 \begin{equation*}
 \Gamma \eqdef \operatorname {Mod} _{-c(\omega _{\Tree})}\sum_{s\in\Tree}\ip  {\mathbf 1}_G, \varphi_s, \varphi_s. 
 \end{equation*}
 This is an assumption that can be assumed by an appropriate modulation of 
 the fixed $ L^2$ function $ f$. 
In the definition of $ \Gamma $, it
is useful to us that we only use the ``smooth'' functions $\varphi_s$
in the definition of this function. 
 Note that $\norm \Gamma .2.{}\lesssim{}\sigma\sqrt{\abs{I_\Tree}}$, which is a consequence of the definition of size
 and the (near) orthogonality of the functions $\varphi_s $ in the case of $ +$tree.

 The purpose of these next remarks is to relate the sums over a $+$tree to oscillation. 
 Recall  that the oscillation is defined relative to a sequence of integers $k_j$. 
For each $J$, consider  $x\in J $ and integers $ \ell$ such that   
$\max(2\abs J,2 ^{\ell_{j-}(x)})<2^\ell<2 ^{\ell _{j+ }(x) }$. 
   We have 
 \begin{equation*}
 \sum_{\substack{s\in\Tree\\ \abs{I_s}=2^\ell }}\ip {\mathbf 1}_G, \varphi_s, f_{s,j}(x)={}
 	\sum_{\substack{s\in\Tree\\  \abs{I_s}=2^\ell}}\ip {\mathbf 1}_G, \varphi_s, \varphi_s(x)\alpha_j(x)
 \end{equation*}
 This is because all of the intervals $\omega_{s+}$ are nested and must contain $\omega_\Tree$, and  if $N(x)\in\omega_{s+}$, then it must also 
 be in every other $\omega_{s'+}$ that is the same size or larger.  What is significant here is that on the right we have a particular 
 scale of (a modulation of) the sum that defines $\Gamma $.

 Furthermore, consider the functions  
 \begin{equation*}
 \Gamma _{j,J}(x)\eqdef{} \operatorname {Mod} _{-c(\omega _{\Tree})}\sum_{\substack{s\in\Tree\\ 
 \max(2\abs J, 2 ^{\ell_{j-}(x)})\le\abs{I_s}\le 2 ^{\ell _{j+ (x)}} } }\ip {\mathbf 1}_G, \varphi_s, \varphi_s. 
\end{equation*}
In particular,
we can choose $\zeta$ as in the definition of our oscillation operator 
(\ref{e.oscWithDenseDef}) so that 
\begin{equation*}
\operatorname {Dil}^{(1)}_{ 2 ^{\ell} } \zeta * \Gamma =
\operatorname {Mod} _{-c(\omega _{\Tree})}
\sum_{\substack{s\in\Tree\\  \abs{ I_s}\ge2 ^{\ell}}}\ip  {\mathbf 1}_G, \varphi_s, \varphi_s .
\end{equation*}
Therefore, we have 
\begin{equation*}
 \Gamma _{j,J}=\left[\operatorname {Dil}^{(1)}_{ \max(2\abs J, 2 ^{\ell_{j-}(x)}) }
 \zeta-\operatorname {Dil}^{(1)}_{ 2 ^{\ell_{j+}(x)}+1 }\zeta\right]*\Gamma .
\end{equation*}

   We conclude that for $x\in E(J)$, 
 \begin{align*} 
\sum_{j=1}^\infty\ABS{ \sum _{\substack{ s\in \Tree \\ 2 \abs J< \abs{I_s} } }  \ip \mathbf 1 _{G}, \varphi_s , f_{s,j}(x)}  & {}\le{}
   \bigg(\sum _{j=1}^\infty \abs{ \Gamma _{j, J}(x) }^2\bigg)^{\frac12}\bigg(\sum_{j=1}^\infty |\alpha_j(x)|^2\bigg)^{\frac12}
   \\&
   {}\lesssim{} \operatorname {Osc}_\delta \Gamma  (x),
 \end{align*} 
 where  we are using the oscillation operator defined in (\ref{e.oscWithDenseDef}).  
 We are able to use this operator here since $2 \abs J< \abs{I_s}$ and $3J$ does not contain any $I_s$, which implies that $J\subset 3I_s$.

 The conclusion of this proof is now at hand. By Lemma~\ref{l.oscWithDense} we have 
 \begin{align*}
 \sum_{\substack{J\in\mathcal J\\ \abs{J}\le3\abs{I_\Tree}}} 
 \int_{E(J)} 
 \sum_{j=1}^\infty\ABS{ \sum _{\substack{ s\in \Tree \\ 2 \abs J< \abs{I_s} } }  \ip \mathbf 1 _{G}, \varphi_s , f_{s,j}(x)} \; dx 
  	{}\lesssim{}& \int_{\bigcup_{\abs{J}\le3\abs{I_\Tree}} E(J)}
	\abs{\operatorname {Osc}_\delta \Gamma(x)}\; dx
 \\{}\lesssim{}& \ABs{\bigcup_{\abs{J}\le3\abs{I_\Tree}} E(J)}^{1/2} 
 \norm\operatorname {Osc}_\delta \Gamma.2.
 \\{}\lesssim{}& \delta\sqrt{\abs{I_\Tree}}\, \norm  \Gamma.2.
 \\{}\lesssim{}& \sigma\delta {\abs{I_\Tree}}. 
\end{align*}

\section{Concluding Remarks} 

Let us pose a conjecture concerning the kernel $ J_H(y) \eqdef \mathbf 1 _{[-1,1]} y ^{-1}$, 
that is the Hilbert transform kernel with a sharp cut off.  

\begin{conjecture}\label{j.osc}  We have the inequality valid for all $ n\ge1$
\begin{equation*}
\norm \operatorname {Osc}_n(J_H; f).p. \lesssim \norm f.p.,\qquad 1<p<\infty .
\end{equation*}
In fact, the implied constant can be taken independent of $ n$. 
\end{conjecture}

The proof as currently presented doesn't permit the deduction of this.  Given the 
central role the Fourier  transform plays in our proof, the technical difficulty 
we come to has a succinct description in terms of $ \widehat {J_H} $.  
Namely, the variation of $ \widehat {J_H} $ is infinite.  But as the variation 
is only logarithmically infinite, one suspects that a proof of the conjecture above 
would have to revisit the proof of Carleson's theorem, with this example in mind.

\bigskip

\begin{corollary}\label{c.discrete} For any measure preserving system
 $ (X,\mu ,\operatorname T)$ and $ f\in L^p (X,\mu)$ for $ 1<p<\infty $, there is a set 
 $ X_f$ of probability one for which for all $ x\in X _{f}$ 
 \begin{equation*}
\lim _{N\to \infty }\sum _{0<\abs k<N} \frac {\operatorname e ^{i \theta k}}k
f (\operatorname T ^{k}x) \quad \text{exists for all $ \theta $.}
\end{equation*}
\end{corollary}

The proof would begin by transferring the oscillation inequality in Theorem~\ref{t.osc}, 
valid on $ \mathbb R $, to the integers $ \mathbb Z $.  This kind of transference 
can be done directly; it is also possible that the necessary result follows 
from known transference results such as Auscher and Carro \cite{MR94b:42007}.  Details 
are left to the reader. 

Likewise, the method of proof that we employ throughout the paper could be adapted to 
shed light on more general singular integrals, as well as the original Wiener--Wintner 
Theorem.  Indeed, an oscillation result could be proved for the latter theorem. 
We do not however pursue these lines here.

\bigskip

The Wiener--Wintner Theorem has a deep extension to the 
Return Time Theorem of Bourgain \cite{MR939436}, see also the appendix to 
\cite{MR1019960}.  This Theorem, which we don't recall in detail here, 
has certain extensions and variants that are currently only approachable 
via  the phase plane methods of the type used in this paper.  The Return 
Time is however a more sophisticated result, and the phase plane methods 
required are correspondingly more difficult.  These issues will be 
explored in forthcoming papers of C.~Demeter, M.~Lacey, T.~Tao, and C.~Thiele.

\end{document}